\documentclass[11pt]{article}
\usepackage{amssymb,amsmath}
\numberwithin{equation}{section}
\date{}
\title{Quantum Double of Yangian of Lie Superalgebra $A(m,n)$ and computation of Universal 
$R$-matrix}
\author{Stukopin V.}
\begin{document}
\maketitle
\vspace{\baselineskip}
\vspace{\baselineskip}
\setcounter{section}{-1}
\setcounter{equation}{0}

\section{Introduction}

Last time along with Yangians of simple (and reductive) Lie algebras it is became studied 
the Yangians of Lie superalgebras of classical type  (see \cite{N}, \cite{N1}, \cite{St}).

  Itself notion of Yangian of simple Lie algebra was introduced by V.Drinfel'd as a 
quantization Lie bialgebra of polynomial currents (with values in this simple Lie algebra) 
and with coalgebra structure defined rational Yang $r$-matrix. 
But dual the Yangian object (for linear algebra Lie $gl(n)$))  was became studied earlier 
in framework Quantum Invers Scaterring Method (QISM).  The V. Drinfel'd shows that this 
object is isomorphic the Yangian. In many papers is ised namely this assignment of Yangian in 
terms of generators which are matrix elements of irredicible representations of Yangian 
in sense of Drinfel'd. As noted above this two approaches essentually isomorphic and its 
employment is dictated of solving problems. In the article  \cite{St} it  was defined Yangian 
of the Lie superalgebra $A(m,n)$ type in the framework Drinfel'd approach and it was 
formulated Poincare-Birkgoff-Witt theorem (PBW-theorem) and theorem on existence 
pseudotriangular structure that is theorem on existence of universal $R$-matrix. This article is natural continuation of the article \cite{St} and its final result is explicite formula 
for universal $R$-matrix. As a corollary we also receive such formula and for the partial 
case of Yangian $Y(sl_2)$ Lie algebra  $sl_2$. I remined that universal $R$-matrix of the 
Yangian $Y(\mathfrak{g})$ of the simple Lie algebra $\mathfrak{g}$ was introduced by 
V.G. Drinfel'd (see \cite{Dr}, \cite{Dr2}) as a formal power series 
 $\it{R}(\lambda)= 1 + \sum_{k=0}^{\infty}\it{R}_k\lambda^{-k-1}$ with coefficients 
$\it{R}_k \in Y(\mathfrak{g})^{\otimes2}$, which conjugates the comultiplication  
$\Delta$ and opposite comultiplication $\Delta' = \tau \circ \Delta, \tau (x \otimes y = 
y \otimes x)$ ($\tau (x \otimes y = (-1)^{deg(x)deg(y)}y \otimes x)$ for Lie superalgebras. 
Explicite defintions will be done later in the text of this article.) More exactly, 
$\it{R}(\lambda)$ conjugate the images $\Delta$ and  $\Delta'$ under action operator 
$id \otimes T_{\lambda}$, where $T_{\lambda}$  is a quantum counterpart of shift operator, аnd  
$id$ be an identical operator. The $\it{R}(\lambda)$ behave as it is image of some hypotetical 
R-matrix $R$ under action $id \otimes T_{\lambda}$ conjugates $\Delta$  and $\Delta'$. 
Drinfel'd called such formal power series $\it{R}(\lambda)$ the pseudotriangular structure and 
proved it existence for $Y(\mathfrak{g})$, when $\mathfrak{g}$ be a simple Lie algebra. 
But explicite formula for $\it{R}(\lambda)$ it  hasn't received up to now. 
If we shall see on classical counterparts of notions  $\it{R}(\lambda)$ and $R$, namely, on  
classical $r$-matrices $\it{r}(\lambda)$ and $r$, then $r$ be an element of a topological tensor square of a classical double and $\it{r}(\lambda)=(id \otimes T_{\lambda})r$, where 
$T_{\lambda}f(u)=f(u+\lambda)$  be a shift operatorа. Then we can naturally to expect that and in the quantum case $\it{R}(\lambda)$ will be an image of universal R-matrix $R$ of quantum double of Yangian  under the action of some shift operator.  
When V.Drinfel'd defined pseudotriangular structure he didn't know good description of 
Yangian double in terms of generators and defining relations and universal R-matrix of Yangian 
double. But in the middle of 90-th S.Khoroshkin and V.Tolstoy receeived the description of 
Yangian double in the terms of generators and defining relations and they computed the 
multiplicative formula of universal R-matrix of Yangian Double (see \cite{Kh-T}).

In this article we describe quantum double  $DY(\mathfrak{g})$ of Yangian of Lie superalgebra 
$\mathfrak{g}=A(m,n)$ in the terms of generators and defining relations. We also calculate the 
universal  R-matrix of Yangian Double follow the plan suggested in article \cite{Kh-T}. 
The main result of this article is a such formula for universal R-matrix for $DY(A(m,n))$. 
This formula is represented in the factorable form as an product of three factors each of them 
is an infinite ordered product. It should be mentioned that computation of universal R-matrix 
of Yangian Double based on the same ideas as a computation of universal R-matrix of quantized 
universal enveloping algebra of affine Lie algebra  (see \cite{L-S-S}, \cite{T-Kh}). 
In the \cite{L-S-S} it was quantum Weyl group for computation the multiplicative formula of 
universal R-matrix. In the case of the Yangian Double we havn't full counterpart quantum Weyl 
group. But partial analogies we use completely. Namely, operator $t^{\infty}$  we can consider 
as an counterpart of longest element of affine Weyl group. We also can interpret the twist 
$F$ which can use for construction of universal R-matrix in the terms of counterpert of elements of affine Weyl group. After them as a formula of universal R-matrix of Yangian double is received we calculate the universal R-matrix of Yangian applying to receiving formula the operator $id \otimes T_{\lambda}$. Further calculation bases on the description action this operators $id \otimes T_{\lambda}$  on the generators of dual Hopf superalgebra to Yangian in quantum  double.

\section{Quantum Yangian Double of Lie Superalgebra $A(m,n)$.}

	I recall that the Yangian $Y(\mathfrak{g})$  of basic Lie superalgebra  $\mathfrak{g}$
(see \cite{F-S}, \cite{K}) is a deformation of universal enveloping superalgebra  $U(\mathfrak{g}[t])$
of bisuperalgebra Lie $\mathfrak{g}[t]$ of polynomial currents. 
The structure of the Lie bisuperalgebra is defined by cocycle 
$\delta: \mathfrak{g} \rightarrow \mathfrak{g} \bigwedge \mathfrak{g}$

\begin{equation}
\delta: a(u) \rightarrow [a(u) \otimes 1 + 1 \otimes a(v),
r(u,v)], \label{coc1}
\end{equation}

where $$r(u,v)= \frac{\mathfrak{t}}{u-v},$$  and  $\mathfrak{t}$  be a Casimir operator, defined   
nondegenerate scalar product $(\cdot, \cdot)$  on basic Lie superalgebra $\mathfrak{g}$ (which 
exist on very basic Lie superalgebra (see \cite{F-S}, \cite{K})). 
Other words, let $\{e_i\}, \{e^i\}$ be dual bases in $\mathfrak{g}$  relatively this scalar product. 
Then $\mathfrak{t}= \sum_i e_i \otimes e^i$. Further, let $\mathfrak{g}=A(m,n)$. 
The Lie superalgebra $\mathfrak{g}$ as an each basic Lie superalgebra is defined by itself Cartan 
matrix $A=(a_{i,j})_{i,j=1}^{m+n+1}$.

Nonzero elements of Cartan matrix are follows: $a_{i,i}=2, a_{i,i+1}=a_{i+1,i}=-1, i < m+1; a_{i-1,i}=a_{i,i-1}=1, a_{i,i}=-2,
m+1<i, i \in I=\{1, \cdots , m+n+1\}.$  Then Lie superalgebra  $\mathfrak{g}$ is generated by the 
generators $h_i, x^{\pm}_i, i \in I,$ where generators $x^{\pm}_{m+1}$ are odd and other generators are 
even. These generators satisfy the following defining relations:

\begin{eqnarray}
&[h_i,h_j]=0
&[h_i,x_j^{\pm}]=\pm a_{ij}x_j^{\pm}, \quad\\
&[x_i^+, x_j^-]=\delta_{ij} h_i, \quad\\
&[[x^{\pm}_{m+1},x^{\pm}_{m+2}], [x^{\pm}_{m+2}, x^{\pm}_{m+1}]]=0\\
&[x_i^{\pm},[x_i^{\pm}, x_j^{\pm}]]=0 \quad
\end{eqnarray}

As usual, $[\cdot, \cdot]$ denotes supercommutator: $[a,b]= ab - (-1)^{p(a)p(b)}ba$.
Let $\Pi = \{\alpha_1, \cdots, \alpha_{m+1}, \cdots, \alpha_{m+n+1}\}$  be a set of 
simple roots, $\Delta (\Delta_+)$  be a set of all roots (positive roots).
Let also  $\{x_{\alpha}, x_{-\alpha}\}, \alpha \in \Delta_+$  be a Cartan-Weyl base,
normalized by condition  $(x_{\alpha}, x_{-\alpha}\})=1$.  We shall use notation 
$(\alpha_i, \alpha_j) = a_{ij}$.

Below, we also use notation $\mathfrak{g}:= A(m,n)$.


 {\bf Definition 1.1.} (see. \cite{St}) Yangian $Y(\mathfrak{g})$ of Lie superalgebra 
$\mathfrak{g}$  be a Hopf superalgebra over  $\cal{C}$, generated as an associative superalgeba 
by generators $h_{i,k}:= h_{\alpha_i, k}, x^{\pm}_{i,k}:=x^{\pm}_{\alpha_i,k}, i \in I, k \in Z_+$, which satisfy  the following defining relations:
\begin{eqnarray}
&[h_{i,k},h_{j,l}]=0, \quad \label{2.1}\\
&\delta_{i,j}   h_{i,k+l}=[x_{i,k}^+,x_{j,l}^{-} ], \label{2.2}\\
&[h_{i,k+1},x_{j,l}^{\pm}]=[h_{i,k},x_{j,l+1}^{\pm}]+
(b_{ij}/2)(h_{i,k}x_{j,l}^{\pm}+x_{j,l}^{\pm} h_{i,k}), \label{2.3}  \\
&[h_{i,0},x_{j,l}^{\pm}] =  \pm b_{ij}x_{j,l}^{\pm}, \label{2.4}\\
&[x_{i,k+1}^{\pm},x_{j,l}^{\pm}]=[x_{i,k}^{\pm},x_{j,l+1}^{\pm}]+
(b_{ij}/2)(x_{i,k}^{\pm}x_{j,l}^{\pm}+x_{j,l}^{\pm}x_{i,k}^{\pm}), \label{2.5}  \quad \\
& \sum_{\sigma}[x_{i,k_{\sigma (1)}}^{\pm}, \cdots
[x_{i,k_{\sigma (r)}}^{\pm},x_{j,l}^{\pm}]...]=0,
i\neq j,  r=n_{ij}=2  \label{2.6}\\
&[[x^{\pm}_{m,k},x^{\pm}_{m+1,k}],[x^{\pm}_{m+2,k},
x^{\pm}_{m+1,k}]]=0 \label{2.7}
\end{eqnarray}

 The sum taken over all permutations $\sigma$ of set $\{1,...,r\}.$  Parity function take the 
following values on generators:
$p(x_{j,k}^{\pm})=0,$ for $k\in Z_+, j \in I\setminus \tau$ $p(h_{i,k})=0,$
for  $i \in I, k\in Z_+,$  $p(x_{i,k}^{\pm})=1, k \in Z_+, i \in \tau.$
  
 Comultiplication on generators $h_{i,k}, x^{\pm}_{i,k}, i \in I, k =0,1$ is defined by 
following formulas: 
\begin{eqnarray}
&\Delta(x) = x \otimes 1 + 1 \otimes x, x \in \mathfrak{g} \label{2.8}   \\
&\Delta(h_{i,1})= h_{i,1} \otimes 1 + 1 \otimes h_{i,1} + [h_{i,0} \otimes 1, \mathfrak{t}_0] +  h_{i,0}\otimes h_{i,0} = \nonumber \\
&h_{i,1} \otimes 1 + 1 \otimes h_{i,1}+ h_{i,0}\otimes h_{i,0} - \sum _{\alpha \in \Delta_+} (-1)^{deg(x_{\alpha})} (\alpha_i, \alpha) x_{-\alpha} \otimes x_{\alpha}; \label{2.9}\\
&\Delta(x^-_{i,1})= x^-_{i,1} \otimes 1 + 1 \otimes x^-_{i,1} + [1 \otimes x^-_{i,0}, \mathfrak{t}_0] =  \nonumber \\
&x^-_{i,1} \otimes 1 + 1 \otimes x^-_{i,1} + \sum _{\alpha \in \Delta_+} (-1)^{deg(x_{\alpha})} [x_{-\alpha_i}, x_{-\alpha}] \otimes x_{\alpha};\label{2.10}\\
& \Delta(x^+_{i,1})= x^+_{i,1} \otimes 1 + 1 \otimes x^+_{i,1} + [x^+_{i,0} \otimes 1, \mathfrak{t}_0] =  \nonumber \\
&x^+_{i,1} \otimes 1 + 1 \otimes x^+_{i,1} - \sum _{\alpha \in \Delta_+} (-1)^{deg(x_{\alpha})} x_{-\alpha} \otimes [x_{\alpha_i},x_{\alpha}];\label{2.11} \end{eqnarray}

   Let's note that universal enveloping superalgebra $U(\mathfrak{g})$ naturally embedded in 
$Y(\mathfrak{g})$.

	Let's introduce the quantum double $DY(\mathfrak{g})$ of Yangian $Y(\mathfrak{g})$. 
I recall the definition of quantum double  (see \cite{Dr1}). Let $A$  be a Hopf superalgebra. 
Let's denote by $A^0$ the dual Hopf superalgebra $A^*$ with opposite comultiplication. Then 
quantum double $DA$ of Hopf superalgebra $A$ be a such quasitriangular Hopf superalgebra 
$(DA, R)$, that $DA$ contains $A, A^0$ as a Hopf subsuperalgebras; $R$ be an image of canonical 
element of $A \otimes A^0,$ corresponding the identical operator under embedding in $DA \otimes DA$; 
linear map $A \otimes A^0 \rightarrow DA, a \otimes b \rightarrow ab$  be a bijection.

Let's note if Hopf superalgebra $A$ be a quantization of bisuperalgebra Lie $\mathfrak{g}$, then 
quantum double $DA$ of $A$ be a quantization of a classical double 
$\mathfrak{g} \oplus \mathfrak{g}^*$ of a Lie bisuperalgebra $\mathfrak{g}$. Moreover the 
cobracket in a classical double is defined by formula: 
$\delta = \delta_{\mathfrak{g}} \oplus (-\delta_{\mathfrak{g}^*}).$

 Let $C(\mathfrak{g})$ (see \cite{Dr3}, \cite{Kh-T})  be an associative superalgebra generated by 
generators $h_{i,k}, x^{\pm}_{i,k}, i \in I, k \in Z$, which satisfy above mentioned defining 
relations (\ref{2.1})-(\ref{2.7}).
 If to define the degrees of generators of $C(\mathfrak{g}$ by formula:
$deg(h_{i,k})=deg(x^{\pm}_{i,k})=k$, then we receive the following filtration on $C(\mathfrak{g})$:
\begin{equation}
\cdots C_{-n} \subset \cdots \subset C_{-1} \subset C_0 \subset \cdots \subset C_m \subset \cdots C(\mathfrak{g}),
\end{equation}
where $C_k= \{ x \in C(\mathfrak{g}): deg(x) \leq k\}$.

Let  $\bar{C}(\mathfrak{g})$ be a formal completion of $C(\mathfrak{g})$ relatively this 
filtration. The generators $x^{\pm}_{i,k}, h_{i,k}, i \in I, k \in Z_+$ generate Hopf subsuperalgebra 
$Y^+(\mathfrak{g})$ в $\bar{C}(\mathfrak{g})$, isomorphic to $Y(\mathfrak{g})$.  Let 
$Y^-(\mathfrak{g})$ be a closed subsuperalgebra in $\bar{C}(\mathfrak{g})$, generated by generators
$x^{\pm}_{i,k}, h_{i,k}, i \in I, k < 0.$

{\bf  Theorem 1.1.}  Hopf Superalgebra $Y^0(\mathfrak{g})$ isomorphic to $Y^-(\mathfrak{g})$.

	This theorm will be follows from the results which will be formulated below. From 
theorem 1.1 it follows that Hopf superalgebra $Y^-(\mathfrak{g})$ be a quantization of 
bisuperalgebra Lie $t^{-1} \mathfrak{g}[[t^{-1}]]$ (with cocycle (\ref{coc1})).

For description $DY(\mathfrak{g})$ it is convenient to introduce the generating functions 
("fields")  $e^+_i(u):= \sum_{k \geq 0} x^+_{i,k}u^{-k-1}, e-_i(u):= -\sum_{k < 0} x^+_{i,k}u^{-k-1},
f^+_i(u):= \sum_{k \geq 0} x^-_{i,k}u^{-k-1}, h^+_i(u):= 1 + \sum_{k \geq 0} h_{i,k}u^{-k-1},
f^-_i(u):= -\sum_{k < 0} x^-_{i,k}u^{-k-1}, h_i(u):= 1 - \sum_{k < 0} h_{i,k}u^{-k-1}.$

{\bf Proposition 1.1.} Defining relations (\ref{2.1})-(\ref{2.7}) in superalgebra $\bar{C}(\mathfrak{g})$
are equivalent the following relations for generating functions
\begin{eqnarray}
&[h^{\pm}_i(u), h^{\pm}_j(v)]=0, [h^+_i(u), h^-_j(u)]=0,\\
&[e^{\pm}_i(u), f^{\pm}_j(v)]=-\delta_{i,j}\frac{h^{\pm}_i(u)- h^{\pm}_i(v)]}{u - v},\\
&[e^{\pm}_i(u), f^{\mp}_j(v)]=-\delta_{i,j}\frac{h^{\mp}_i(u)- h^{\pm}_i(v)]}{u - v},\\
&[h^{\pm}_i(u), e^{\pm}_j(v)]=-\frac{(\alpha_i, \alpha_j)}{2} \frac{\{h^{\pm}_i(u), (e^{\pm}_j(u)- e^{\pm}_j(v))\}}{u - v}, \quad\\
&[h^{\pm}_i(u), e^{\mp}_j(v)]=-\frac{(\alpha_i, \alpha_j)}{2} \frac{\{h^{\pm}_i(u), (e^{\pm}_j(u)- e^{\mp}_j(v))\}}{u - v},\\
&[h^{\pm}_i(u), f^{\pm}_j(v)]=\frac{(\alpha_i, \alpha_j)}{2} \frac{\{h^{\pm}_i(u), (f^{\pm}_j(u)- f^{\pm}_j(v))\}}{u - v},\quad\\
&[h^{\pm}_i(u), f^{\mp}_j(v)]=\frac{(\alpha_i, \alpha_j)}{2} \frac{\{h^{\pm}_i(u), (e^{\pm}_j(u)- e^{\mp}_j(v))\}}{u - v},\quad\\
&[e^{\pm}_i(u), e^{\pm}_j(v)] + [e^{\pm}_j(u), e^{\pm}_i(v)] =-\frac{(\alpha_i, \alpha_j)}{2} \frac{\{(e^{\pm}_i(u)- e^{\pm}_i(v)), (e^{\pm}_j(u)- e^{\pm}_j(v))\}}{u - v}, \qquad \\
&[e^+_i(u), e^-_j(v)] + [e^+_j(u), e^-_i(v)] =-\frac{(\alpha_i, \alpha_j)}{2} \frac{\{(e^+_i(u)- e^-_i(v)), (e^+_j(u)- e^-_j(v))\}}{u - v}, \qquad\\
&[f^{\pm}_i(u), f^{\pm}_j(v)] + [f^{\pm}_j(u), f^{\pm}_i(v)] =-\frac{(\alpha_i, \alpha_j)}{2} \frac{\{(f^{\pm}_i(u)- f^{\pm}_i(v)), (f^{\pm}_j(u)- f^{\pm}_j(v))\}}{u - v},\quad\\
&[f^+_i(u), f^-_j(v)] + [f^+_j(u), f^-_i(v)] =-\frac{(\alpha_i, \alpha_j)}{2} \frac{\{(f^+_i(u)- f^-_i(v)), (f^+_j(u)- f^-_j(v))\}}{u - v},\qquad\\
&[e_i^{\epsilon_1}(u_1), [e_i^{\epsilon_2}(u_2), e_j^{\epsilon_3}(u_3)]] + [e_i^{\epsilon_2}(u_2), [e_i^{\epsilon_1}(u_1), e_j^{\epsilon_3}(u_3)]] =0 ,\qquad \\
&[f_i^{\epsilon_1}(u_1), [f_i^{\epsilon_2}(u_2), f_j^{\epsilon_3}(u_3)]] + [f_i^{\epsilon_2}(u_2), [f_i^{\epsilon_1}(u_1), f_j^{\epsilon_3}(u_3)]] =0 ,\qquad \\
&[[e_m^{\epsilon_1}(u_1), e_{m+1}^{\epsilon_2}(u_2)], [e_{m+2}^{\epsilon_3}(u_3), e_{m+1}^{\epsilon_4}(u_4)]]=0,\quad \\
&[[f_m^{\epsilon_1}(u_1), f_{m+1}^{\epsilon_2}(u_2)], [f_{m+2}^{\epsilon_3}(u_3), f_{m+1}^{\epsilon_4}(u_4)]]=0.\quad
\end{eqnarray}

\section{Triangular decomposition and pairing formulas.}

        Let $Y'_+, Y'_0, Y'_-$  be subsuperalgebras  (without unit) in $Y(\mathfrak{g})$,
generating of elements $x^+_{ik}, h_{ik}, x^-_{ik}, (i \in I, k \in Z_+),$
correspondingly.
Let $Y_+, Y_0, Y_-$ be subsuperalgebras $Y'_+, Y'_0, Y'_-$ with unit element.

{\bf  Proposition 2.1.} Multiplication in $Y(\mathfrak{g})$ induces isomorphism of vector 
superspaces 

\begin{equation}
Y_+ \otimes Y_0 \otimes Y_- \rightarrow Y(\mathfrak{g})
\end{equation}

  This proposition is partial case of theorem 3 from \cite{St}.  Let's  extend this proposition 
on  $DY(\mathfrak{g})$.

For these we need some simple properties of comultiplication on  $Y(\mathfrak{g})$, which is proved 
by induction using formulas (\ref{2.8})-(\ref{2.11}) and relations (\ref{2.1})-(\ref{2.7}), and using also 
that fact that comultiplication be homomorphism of associative superalgebra, i.e.  
$\Delta(a \cdot b)= \Delta(a) \cdot \Delta(b)$.

{\bf  Proposition 2.2.}
1) $\Delta(x) = x \otimes 1 (mod Y \otimes Y'_+)$, for all $x \in Y'_+;$\\
2) $\Delta(y) = 1 \otimes y (mod Y'_- \otimes Y)$, for all $y \in Y'_-.$\\

{\bf Corollary.} 1) $\Delta(Y_+) \subset Y\otimes Y_+;$\\
2) $\Delta(Y_-) \subset Y_-\otimes Y.$

  So, we  have that $Y_+ (Y_-)$ be a right (left) coideal in $Y = Y(\mathfrak{g})$.

     Let's also $BY'_{\pm}$ be a subsuperalgebra (without unit) in $Y(\mathfrak{g})$,
generated by $x^{\pm}_{ik}, h_{jr}, (i, j \in I, k,r \in Z_+).$

{\bf  Proposition 2.3.}
1) $\Delta(e) = e \otimes 1 (mod Y \otimes BY'_+)$, for all $e \in BY'_+;$\\
2) $\Delta(f) = 1 \otimes f (mod BY'_- \otimes Y)$, for all $f \in BY'_-.$\\
3) $\Delta(h) = h \otimes 1 (mod Y \otimes BY'_+)= 1 \otimes h (mod BY'_- \otimes Y)$,
for all  $f \in Y'_0.$\\

    Properties  1), 2) is proved also as analogous properties in  proposition 2.2, property 3) 
follows from 1), 2).

    Let $<\cdot, \cdot> : Y(\mathfrak{g}) \otimes Y^0(\mathfrak{g}) \rightarrow C$  be a 
canonical bilinear pairing $Y(\mathfrak{g})$  and its dual Hopf superalgebra $Y^*(\mathfrak{g})$ 
with opposite comultiplication.
(We denote by $Y^0(\mathfrak{g})$ the $Y^*(\mathfrak{g})$ with opposite comultiplication.)  
From definition imply the next properties of this pairing.

$<xy, x'y'>= <\Delta(xy), x' \otimes y'> =
(-1)^{p(x)p(y)}<y \otimes x, \Delta(x'y')>, \\
<x \otimes y><x^{\prime}\otimes y^{\prime}> =
(-1)^{p(x)p(y)}<x, x^{\prime}><y, y^{\prime}>,$ \\
для
$\forall x, y \in  Y(\mathfrak{g}), \forall x^{\prime}, y^{\prime} \in
Y^0(\mathfrak{g})$

    Let $A, B$  are subsuperalgebras of $Y(\mathfrak{g})$. Let's also 
$(AB)_{\bot} := \{x' \in Y^0(\mathfrak{g}): <ab, x'>=0,
\forall a \in A, b \in B\}$.
It is easy to check that 
$(Y \cdot BY'_-)_{\bot}, (BY'_+ \cdot Y)_{\bot}, (Y \cdot Y'_-)_{\bot},
(Y'_+ \cdot Y)_{\bot}$  are subsuperalgebras of $Y^0(\mathfrak{g})$.  Let \\
$Y^*_+ := (Y \cdot BY'_-)_{\bot}, BY^*_+ := (Y \cdot Y'_-)_{\bot},
Y^*_- := (BY'_+ \cdot Y)_{\bot}, \\
(BY)^*_- := (Y'_+ \cdot Y)_{\bot}, Y^*_0 := BY^*_+ \bigcap BY^*_-.$

    {\bf  Proposition 2.4.} 1) For all  $x \in Y_+, h \in Y_0, y \in Y_-, x' \in Y^*_+, h' \in Y^*_0, y' \in Y^*_-$
canonical pairing is factorized as 
$$<xhy, x'h'y'> = (-1)^{deg(x')deg(y)}<x, x'><h, h'><y, y'>.$$
2) Multiplication in $Y^0(\mathfrak{g})$ induces isomorphism of vector spaces:\\
$Y^*_+ \otimes Y^*_0 \otimes Y^*_- \rightarrow Y^0(\mathfrak{g}).$\\
3) PBW-theorem is fulfilled for $Y^0(\mathfrak{g})$.

{\bf  Proof.} Let' prove 1). \\
$<xhy, x'h'y'>= <\Delta(xh)\cdot \Delta(y), x'h' \otimes y'>= <\Delta(x) \Delta(h) \Delta(y), x'h' \otimes y'> =\\
<(x \otimes 1 + \sum a_n \otimes x_n)(h \otimes 1 + \sum \tilde{a}_s \otimes \tilde{x}_s) (1 \otimes y + \sum y_m \otimes a'_m), x'h' \otimes y'>= \\
<xh \otimes y, x'h' \otimes y'> + < \sum c_r \otimes d_r, x'h' \otimes y'> = \\
(-1)^{deg(x')deg(y)},xh, x'h'> <y, y'>$.

Let's noticed that $<\sum c_r \otimes d_r, x'h' \otimes y'> = \sum_r (-1)^{deg(x')deg(d_r)} <c_r, x'h'><d_r, y'>=0.$
As $<d_r, y'>=0$,   $d_r \in Y'_+Y, y' \in (BY'_+Y)_{\bot}$ and we have $<xh, x'h'>= 
< (x \otimes 1 + \sum a_n \otimes x_n)(1\otimes h + \sum y_m \otimes b_m), x' \otimes h'> =
<x \otimes h, x' \otimes h'> + 0 = <x, x'><h, h'>$ therefore we received proposition of 1).

    Let's note that 2) follows from 3). Let' prove 3). Let's choose  PBW base 
$Y(\mathfrak{g})$. Every vector of this base can be represented in the followin form:
$xhy$, where $x \in Y_+, h \in H, y \in Y_-.$  Then biorthogonal vector in view of 1),
can be represented in the form: $x'h'y'$, where $ x' \in Y^*_+, h \in H^*, y \in Y^*_-.$
These vectors also form base in $Y^0(\mathfrak{g})$. This fact proves  3).

    Let's study this pairing in detail. First, let's describe the PBW base for  $Y(\mathfrak{g})$ 
in detail (see alternative description in \cite{St}). 
Let's as above $\Delta, \Delta_+$ denote the set of roots, set of positive roots of Lie superalgebra 
$A(m,n)$.  Let's consider also the set $\hat{\Delta}^{re}$ of real roots of affine (nontwisted)
Lie superalgebra  $A(m,n)^{(1)}$ (see \cite{F-S}). For generators of $DY(\mathfrak{g})$ $x^{\pm}_{i,k}$ 
we shall use the next notation:\\
$x_{\alpha_i + k\delta} := x^{+}_{i,k}, \\
x_{-\alpha_i + k\delta} := x^{-}_{i,k}, i \in I, k \in Z, \alpha_i \in \Delta_+.$

    In this case $\pm\alpha_i + k\delta \in \hat{\Delta}^{re}.$   Let $\Xi \subset \hat{\Delta}^{re}$.
Total linear order $\precneqq$  на $ \Xi$ is called convex (normal), if for all roots 
$\alpha, \beta, \gamma \in \Xi$ such that $\gamma= \alpha + \beta$ we have :\\
$\alpha \precneqq \gamma \precneqq \beta$ or $\beta \precneqq \gamma \precneqq \alpha$.\\
Let's introduce subsets $\Xi_+, \Xi_-$ of set $\hat{\Delta}^{re}$: \\
$\Xi_{\pm} := \{\pm \gamma +k\delta: \gamma \in \hat{\Delta}^{re}_+ \}$.\\

    Let's introduce on $\Xi_+, \Xi_-$ convex orderings $\precneqq_+, \precneqq_-$,
saisfying the following conditions:
\begin{equation}
\gamma + k\delta \precneqq_+ \gamma + l\delta \qquad \mbox{and}\quad -\gamma +
l\delta \precneqq_-  -\gamma + k\delta,
\qquad \mbox{if}\quad k \le l, \quad \mbox{for} \quad \forall \gamma \in \Delta_+  \label{no}
\end{equation}

    Let's define root vectors $x_{\pm \beta}, \beta \in \Xi_+ \cup \Xi_-$ by induction in the 
following way. Let vectors  $x_{\beta_1}, x_{\beta_2}$ are already being constructed. If root 
$x_{\beta_3}$ satisfy conditions:
$x_{\beta_1} \precneqq x_{\beta_3} \precneqq x_{\beta_2}$ and in the segment $(x_{\beta_1}, x_{\beta_2})$ 
we havn't root vecors (which was already being constructed), the let's define root vectors 
$x_{\pm\beta_3}$ by formulas:\\
$x_{\beta_3} = [x_{\beta_1}, x_{\beta_2}], x_{-\beta_3} = [x_{-\beta_2}, x_{-\beta_1}].$

Let's note that convex (normal) ordering connects with natural ordering of elements of 
affine Weyl group in the case of Lie algebras. 

We need the following description of $Y(\mathfrak{g})$, which is an analog of description of quantized 
universal enveloping superalgebra of affine Lie superalgebra. First, let's fix the following convex 
ordering on $\mathfrak{g}=A(m,n)$:\\
$(\epsilon_1 - \epsilon_2, \epsilon_1 - \epsilon_3, \cdots, \epsilon_1 - \epsilon_{m+n+2}),
(\epsilon_2 - \epsilon_3, \epsilon_2 - \epsilon_4, \cdots, \epsilon_2 - \epsilon_{m+n+2}), \cdots
(\epsilon_{m+n+1} - \epsilon_{m+n+2}).$ \\
Here $\epsilon_i - \epsilon_j = \alpha_i + \cdots + \alpha_{j-1}$.

Let's add affine root $\alpha_0$ to the set of simple roots. I remind that  
$\alpha_0 = \delta - \theta, \theta := \alpha_1 + \cdots + \alpha_{m+n+1} = \epsilon_1 - \epsilon_{m+n+1}$ 
be a highest root, $\delta$ be a minimal imaginary root.
Let's consider the next convex ordering on the set $\hat{\Delta}^{re}$ affine real roots:\\
$(\alpha_1, \alpha_1 + \delta, \alpha_1 + 2\delta, \cdots, \alpha_1 + n\delta,
\cdots ),
(\cdots \alpha_1 + \alpha_2 + n \delta,
\cdots, \alpha_1 + \alpha_2 + \delta,  \alpha_1 + \alpha_2),
(\cdots \alpha_1 + \alpha_2+ \alpha_3 + n \delta,
\cdots, \alpha_1 + \alpha_2 + \alpha_3 + \delta, \alpha_1 + \alpha_2 +\alpha_3),
\cdots , (\cdots \alpha_{n+m+1} + n \delta, \cdots, \alpha_{n+m+1} + \delta,
\alpha_{n+m+1})$.

    Let's calculate pairing for root vectors. Let $h^*_{i,k}, e^*_{i,k}, f^*_{i,k}$ are generators 
of $Y^*=Y_-$. Let $e_{i,k} := x^+_{i,k}, f_{i,k} := x^-_{i,k}$.

    {\bf  Proposition 2.5.} The following two conditions are equivalent.\\
1) $<e_{i,k}, e^*_{j, -l-1}> = -\delta_{i,j} \delta_{k,l};$\\
  $<f_{i,k}, f^*_{j, -l-1}> = -\delta_{i,j} \delta_{k,l};$\\
  $ <h_{i,k}, h^*_{j, -l-1}> = -(-\frac{a_{ij}}{2})^{k-l}\frac{a_{ij}k!}{l!(k-l)!} \qquad{for} k \geq l \geq 0.$\\

2)
\begin{eqnarray}
&[h^*_{i,-k}, h^*_{j,-l}] = 0, \quad \\
&\delta_{i,j}   h^*_{i,-k-l} = [e_{i,-k}, f_{j,-l}], \label{equation61}\\
&[h^*_{i,-k-1}, e^*_{j,-l}]=[h^*_{i,-k}, e^*_{j,-l-1}]+
(b_{ij}/2)(h^*_{i,-k}e_{j,-l} + e^*_{j,-l} h^*_{i,-k}),  \\
&[h^*_{i,-k-1}, f^*_{j,-l}] = [h^*_{i,-k}, f^*_{j,-l-1}]-
(b_{ij}/2)(h^*_{i,-k}f_{j,-l}+f^*_{j,-l}h^*_{i,-k}),  \\
&[h^*_{i,0},e^*_{j,l}] =   b_{ij}e^*_{j,l}, \label{equation63}\\
&[h^*_{i,0}, f^*_{j,l}] =  -b_{ij}f^*_{j,l}
\end{eqnarray}
\begin{eqnarray}
&[e^*_{i,-k+1}, e^*_{j,-l}]=[e^*_{i,-k},e^*_{j,-l+1}]+
(b_{ij}/2)\{e^*_{i,-k}, e^*_{j,-l}\},  \quad \\
&[f^*_{i,-k+1},f^*_{j,-l}]=[f^*_{i,-k},f^*_{j,-l+1}]-
(b_{ij}/2)\{f^*_{i,-k}, f^*_{j,-l}\},  \quad \\
& \sum_{\sigma}[e^*_{i,-k_{\sigma (1)}}, \cdots
[e^*_{i,-k_{\sigma (r)}},e^*_{j,-l}]...]=0,
i\neq j,  r=n_{ij}=2\\
& \sum_{\sigma}[f^*_{i,-k_{\sigma (1)}}, \cdots
[f^*_{i,-k_{\sigma (r)}},f^*_{j,-l}]...]=0,
i\neq j,  r=n_{ij}=2\\
&[[e^*_{m,-k_1},e^*_{m+1,-k_2}],[e^*_{m+2,-k_3}, e^*_{m+1,-k_4}]]=0,\\
&[[f^*_{m,-k_1},f^*_{m+1,-k_2}],[f^*_{m+2,-k_3}, f^*_{m+1,-k_4}]]=0.
\end{eqnarray}

{\bf  Proof.}  The proof of this proposition is inconveniently and we having marked basic 
points of proof, omiting technical details. We shall lead the proof by induction on values of 
indexes $k,l$. First of all it is easy to prove next formulas. \\
$\Delta(e_{i,k}) = e_{i,k} \otimes 1 + 1 \otimes e_{i,k} + \sum_{r=0}^{k-1} h_{i,r} \otimes e_{i,k-r} (mod Y Y_- \otimes Y'_+);$ \\
$\Delta(f_{i,k}) = f_{i,k} \otimes 1 + 1 \otimes f_{i,k} + \sum_{r=0}^{k-1} h_{i,r} \otimes e_{i,k-r} (mod Y'_- \otimes Y'_+Y); \\
\Delta(h_{i,k}) = h_{i,k} \otimes 1 + 1 \otimes h_{i,k} + \sum_{r=0}^{k-1} h_{i,r} \otimes h_{i,k-r} (mod Y Y'_- \otimes Y'_+Y);$ \\

From this formulas it is follows the next equalities:\\
$\Delta(e_{i,k}e_{j,l}) = e_{i,k}e_{j,l} \otimes 1 + 1 \otimes e_{i,k}e_{j,l} + e_{i,k} \otimes e_{j,l} +
(-1)^{deg(e_{i,k})deg(e_{j,l})} e_{j,l} \otimes e_{i,k} (mod YY'_- \otimes Y'_+); \\
\Delta(f_{i,k}f_{j,l}) = f_{i,k}f_{j,l} \otimes 1 + 1 \otimes f_{i,k}f_{j,l} + f_{i,k} \otimes f_{j,l} +
(-1)^{deg(f_{i,k})deg(f_{j,l})} f_{j,l} \otimes f_{i,k} (mod Y'_- \otimes Y'_+Y); \\
\Delta(h_{i,k}h_{j,l}) = h_{i,k}h_{j,l} \otimes 1 + 1 \otimes h_{i,k}h_{j,l} + h_{i,k} \otimes h_{j,l} +
h_{j,l} \otimes h_{i,k} (mod YY'_- \otimes Y'_+Y).$ \\
    
      Using these formulas and definition of quantum double we can prove by induction the 
invariance of this pairing on generators of Yangian Double. 

\begin{equation}
<[a,b],c> = <a,[b,c]>  \label{p1}
\end{equation}

	We omit the proof of this fact realizing that proof of such simple fundamental fact 
it must be short and idea's. We have only proof bases on induction using above written formulas 
and next definition of Hopf pairing.    
\begin{eqnarray}
&<ab,cd> = <\Delta(ab),c\otimes d> = (-1)^{deg(a)deg(b)}<b\otimes a, \Delta(cd)>
 \label{s1}\\
& <a, 1>=\epsilon(a)>, <1, b>= \epsilon(b)  \qquad  \label{s2}
\end{eqnarray}

	Now we can to show how the condition 1) follows from condition 2). Let's show, 
for example, how by induction it is  derived pairing formula on Cartan generators of 
Yangian Double from commutative relations using formula  \ref{p1}. For  $m=n=0$, 
proving formulas coincide with its quasiclassical limits for which they are evidently correct. 
Let these formulas correct for $m \ge k, n<l+1$.
Let's show that they correct for $n=l+1$.\\
$<h_{i,k}, h_{j,l}> = -<e_{i,0},[f_{i,k}, h_{j,l}]> = <e_{i,0},
[h_{j,0},f_{i,k-l-1}] + \\
\frac{1}{2}a_{ij}\sum_{s=0}^l\{h_{j,s-l-1}, f_{k-s-1}\}>  =
-\frac{1}{2}a_{ij}(<e_{i,0}, \sum_{s=0}^l ([h_{j,s-l-1}, f_{k-s-1}]
+ \\
2f_{i,k-s-1}h_{j,s-l-1})>) =
-\frac{1}{2}a_{ij}(<e_{i,0}, \sum_{s=0}^l [h_{j,s-l-1},
f_{k-s-1}]>\\
+ 2<e_{i,0}, \sum_{s=0}^l [f_{i,k-s-1}, h_{j,s-l-1}]>)$. \\
   Second summand equal zero in view of inductive assumption. Let's transform first summand. 
Let's decrease degree of right-hand side in pairing formula using defining relations in 
Yangian Double. 

$<h_{i,k}, h_{j,l}> =a_{ij}<e_{i,0}, \sum_{s=0}^l(l+s-1) [h_{j,-s}, f_{k-l+s-2}]a_{ij}/2> =
-(\frac{1}{2}a_{ij})^2<e_{i,0}, \sum_{s=0}^l(l+s-1) [h_{j,-s}, f_{k-l+s-2}]> = \cdots \\
=-(\frac{1}{2}a_{ij})^{k-l}(<e_{i,0},[h_{j,0}, f_{i,-1}]>C^{k-l-1}_{k-l-1} + C^{k-l-1}_{k-l} +
\cdots \\
C^{k-l-1}_{k-l+l-1})=
-(\frac{1}{2}a_{ij})^{k-l}C^{k-l}_{k}a_{ij}.$\\
  First pairing formula is proved. The second formula is proved simpler by analogous arguments.  
	The proving of sufficiency rather inconvin inconveniently and we omit it here. Note, only, 
that actually it is also leaded by induction and  based on formulas (\ref{s1}), (\ref{s2}).

    {\bf  Theorem 2.1.} 1) Subsuperalgebras $Y^*_+, H^*, Y^*_-$ of superalgebra $Y_-$ are 
generated by fields 
$e^-_i(u),  h^-_i(u), f^-_i(u);$\\
2) Pairing of generators of subsuperalgebras $Y_+, Y_-$ of superalgebra $DY(\mathfrak{g})$ is assigned 
the next relations for 
$|v|< 1 <|u|$: \\
\begin{eqnarray}
&<e^+_i(u), f^-_j(v)> = <f^+_i(u), e^-_j(v)> = \frac{\delta_{i,j}}{u-v};     \\
&<h^+_i(u), h^-_j(v)> =  \frac{u-v + \frac{1}{2}(\alpha_i, \alpha_j)}{u-v - \frac{1}{2}(\alpha_i, \alpha_j)}\label{2.2}
\end{eqnarray}

{\bf  Proof.} Theorem imply from proposition 2.5.

\section{Computation of Universal $R$- matrix of Yangian Double $DY(\mathfrak{g})$.}

	First, I remind the definition of universal $R$- matrix for quasitriangular Hopf 
topological superalgebra, which is natural generalization of notion of universal $R$-matrix for 
quasitriangular Hopf algebra  (see \cite{Dr1}).

	Universal $R$-matrix for quasitriangular topological Hopf superalgebra $A$ is called such 
invertible element $R$ from some extension of completion of tensor square $A \hat{\otimes} A$ and satisfying next conditions :\\
$$\Delta^{op}(x) = R \Delta(x) R^{-1}, \qquad \forall x \in A;$$
$$(\Delta \otimes id) R = R^{13}R^{23}, (id \otimes \Delta) R = R^{13} R^{12}, $$
где $\Delta^{op}= \sigma \circ \Delta, \sigma(x \otimes y) = (-1)^{p(x)p(y)} y \otimes x.$

	If $A$ be a quantum double of Hopf superalgebra  $A^+$, i.e.   
$A \cong A^+ \otimes A^-$, $A^-:=A^0$ be a dual to $A$ Hopf superalgebra with opposite 
comultiplication, then  $A$ be a quasitriangular Hopf superalgebra and universal $R$- matrix in 
$A$ assume  next canonical presentation: $$ R = \sum e_i \otimes e^i,$$
where $\{e_i\}, \{e^i\}$ are dual bases in $A^+, A^-$, respectively.

   Let $Y^{\pm}_+, Y^{\pm}_0, Y^{\pm}_-$  are subsuperalgebras in $DY(\mathfrak{g})$,
generated by fields $e^{\pm}_i(u), h^{\pm}_i(u), \\
f^{\pm}_i(u), i \in I,$ respectively.

  {\bf  Proposition 3.1.} 1) Universal $R$-matrix of Yangian Double can be presented 
in the next factorizable form:
$$R = R_+ R_0 R_-,$$
where $R_+ \in Y^+_+ \otimes Y^-_-, R_0 \in Y^+_0 \otimes Y^-_0, R_- \in Y^+_- \otimes Y^-_+.$\\
2)Pairing on the base elements can be computed according to the next formulas:
\begin{eqnarray}
&<e^{n_0}_{\beta_0}e^{n_1}_{\beta_1} \cdots e^{n_k}_{\beta_k}, e^{m_0}_{-\beta_0- \delta}e^{m_1}_{-\beta_1- \delta}
\cdots e^{m_k}_{-\beta_k- \delta}> = \quad \nonumber\\
&(-1)^{n_0 + \cdots + n_k} \delta_{n_0, m_0} \cdots
\delta_{n_0, m_0}\cdot n_0!n_1! \cdots n_k! \cdot \alpha(\gamma_0)^{n_0} \cdots
\alpha(\gamma_k)^{n_k}(-1)^{\theta(\beta_0)+ \cdots + \theta(\beta_k)};\qquad\\
&< e^{n_k}_{-\beta_k} \cdots  e^{n_1}_{-\beta_1}  e^{n_0}_{\beta_0}, e^{m_k}_{-\beta_k- \delta}
\cdots e^{m_1}_{-\beta_1- \delta} e^{m_0}_{-\beta_0- \delta}> = \nonumber\\
&(-1)^{n_0 + \cdots + n_k} \delta_{n_0, m_0} \cdots
\delta_{n_0, m_0}\cdot n_0!n_1! \cdots n_k! \cdot \alpha(\gamma_0)^{n_0} \cdots
\alpha(\gamma_k)^{n_k} (-1)^{\theta(\beta_0)+ \cdots + \theta(\beta_k)}; \qquad
\end{eqnarray}

    Here $\beta_k= \beta'_k + n'_k\delta$, and coefficients $\alpha(\beta)$ can be calculated 
from condition $[e_{\beta}, e_{-\beta}]=\alpha(\beta)h_{\beta'}$.

    From proposition 3.1 follows 

  {\bf  Lemma 3.1.} The elements $R_+,  R_-$ in decomposition of universal $R$-matrix of 
$DY(\mathfrak{g})$ can be presented in the following form 
\begin{eqnarray}
R_+ = & \overrightarrow{\prod}_{\beta \in \Xi_+} \exp(-(-1)^{\theta(\beta)}a(\beta) e_{\beta} \otimes e_{-\beta}), \label{3.3}\\
R_- = & \overleftarrow{\prod}_{\beta \in \Xi_-}\exp(-(-1)^{\theta(\beta)}a(\beta) e_{\beta} \otimes e_{-\beta}),  \label{3.4}
\end{eqnarray}
where product taken according to normal orderings $\precneqq_+, \precneqq_-$, satisfying 
conditions \ref{no}.

Normalizing constants $a(\beta)$ can be found from the following condition:
\begin{eqnarray*}
&[e_{\beta}, e_{-\beta}] = (a(\beta))^{-1} h_{\gamma} \quad \mbox{if}  \quad \beta = \gamma + n\delta \in \Xi_+, \gamma \in \Delta_+(\mathfrak{g}),\quad\\
&[e_{\beta}, e_{-\beta}] = (a(\beta))^{-1} h_{\gamma} \quad \mbox{if}  \quad \beta = \gamma + n\delta \in \Xi_+, \gamma \in \Delta_+(\mathfrak{g}),\quad\\
\end{eqnarray*}

and $\theta(\beta)= deg(e_{\beta})=deg(e_{-\beta})$  denotes parity of element $e_{\pm \beta}$.

     For description of term $R_0$ we need some auxiliary notions.
First of all, let's introduce "logarithmic" generators $\phi^{\pm}_i(u), i=1, \cdots,r$ by formulas 
\begin{equation}
\phi^{+}_i(u):= \sum_{k=0}^{\infty}\phi_{i,k}u^{-k-1} = \ln h^+_i(u);
\phi^{-}_i(u):= \sum_{k=0}^{\infty}\phi_{i,-k-1}u^k = \ln h^-_i(u)
\end{equation}
  
Let's introduce vector-functions 

${\phi}^{\pm}(u)=
\begin{pmatrix}
\phi^{\pm}_1(u)\\ \phi^{\pm}_2(u)\\ \cdots \\
\phi^{\pm}_r(u)
\end{pmatrix}
h^{\pm}(u)=
\begin{pmatrix}
h^{\pm}_1(u)\\ h^{\pm}_2(u)\\ \cdots \\
h^{\pm}_r(u)
\end{pmatrix}$

   From theorem  2.1  implyes pairing formula in the terms of generating vector-functions
\begin{equation}
<((h^+(u))^T, h^-(v))>=
(\frac{u-v + \frac{1}{2}(\alpha_i, \alpha_j)}{u-v -\frac{1}{2}(\alpha_i, \alpha_j)})_{i,j=1}^r
\end{equation}
    Therefore, for generating functions 
$\phi^+_i(u), \phi^-_j(u)$  pairing formula has the following form 
\begin{equation}
<\phi^+_i (u), \phi^-_j (v)> = ln(\frac{u-v + \frac{1}{2}(\alpha_i,\alpha_j)}
{u-v - \frac{1}{2}(\alpha_i,\alpha_j)})  \label{3.5}
\end{equation}

These formulas we can rewritten in the matrix form as 
\begin{equation}
<(\phi^+(u))^T, \phi^-(v)> =
(\ln(\frac{u-v + \frac{1}{2}(\alpha_i, \alpha_j)}
{u-v - \frac{1}{2}(\alpha_i, \alpha_j)}))_{i,j=1}^r
\end{equation}

\begin{equation}
<(\phi^+(u))^T, \phi^-(v)> = (\ln(\frac{u-v + \frac{1}{2}(\alpha_i - \alpha_j)}
{u-v - \frac{1}{2}(\alpha_i - \alpha_j)}))_{i,j=1}^r   \label{3.6}
\end{equation}

     Further calculation we shall conduct  on the scheme suggested in \cite{Kh-T}.
Using this way we can attach to this calculations some gemetrical sense.

     Along with Yangian Double $DY(\mathfrak{g})$ let's consider Hopf superalgebra 
$\widehat{DY}(\mathfrak{g})$, isomorphic to as associative superalgebra to $DY(\mathfrak{g})$, 
but with another comultiplication defined next formulas:
\begin{eqnarray}
&\tilde{\Delta}(h^{\pm}_i(u))=h^{\pm}_i(u) \otimes h^{\pm}_i(u) \quad \label{3.k} \\
&\tilde{\Delta}(e_i(u))=e_i(u) \otimes 1 + h^{-}_i(u) \otimes e_i(u) \quad \\
&\tilde{\Delta}(f_i(u)) =1 \otimes f_i(u) + f_i(u)\otimes h^{+}_i(u) \quad
\end{eqnarray}

Here \\
$e_i(u):= e^+_i(u) - e^-_i(u) = \sum _{k \in Z} e_{i,k} u^{-k-1},\\
 f_i(u):= f^+_i(u) - f^-_i(u) = \sum _{k \in Z} f_{i,k} u^{-k-1}$

   Such comultiplication it was introduced by V.Drinfel'd (\cite{Dr3}) in the case of Yangians 
(and Yangian Doubles) of simple Lie algebras and it convenient by that the pairing formulas 
relatively this comultiplication has a simple form. 
It is possible to check  (see \cite{Kh-T}), that comultiplications  $\Delta$ и $\tilde{\Delta}$
conjugated by limit operator 
$\hat{t}^{\infty}:= lim_{n \rightarrow \infty} \hat{t}^n$,
$\hat{t}(e_{i,k})=e_{k+1}, \hat{t}(f_{i,k})=f_{k-1}, \hat{t}(h_{i,k})=h_{k}$.
In other words,
\begin{equation}
\tilde{\Delta}(x) = lim_{n \rightarrow \infty}(\hat{t}^n \otimes \hat{t}^n)\Delta(\hat{t}^{-n}(x)),
\end{equation}
for $\forall x \in DY(\mathfrak{g})$. (The convergence it is implied in suitable topology of 
$DY(\mathfrak{g})\otimes DY(\mathfrak{g})$.)

     Let $\widehat{DY}^+(\mathfrak{g})$ ($\widehat{DY}^-(\mathfrak{g})$) be a Hopf subsuperalgebra 
of Hopf superalgebra $\widehat{DY}(\mathfrak{g})$, generated by elements 
$e_{i,k}, k \in Z, h_{i,m}, m \in Z_+$ ($f_{i,k}, k \in Z, h_{i,m}, m <0$). Then 
$\widehat{DY}^-(\mathfrak{g})$ isomorphic to dual Hopf superalgebra $(\widehat{DY}^-(\mathfrak{g}))^*$.
From comultiplication formula  (\ref{3.k}) imply that elements $\phi^{\pm}_{i,k}$ are primitive 
elements in $\widehat{DY}(\mathfrak{g})$.

    Let $\Phi^+ = <\phi^{+}_{i,k}: i \in I= \{1, \cdots, r\}, k \in Z_+>,
\Phi^- = <\phi^{-}_{i,-k-1}: i \in I= \{1, \cdots, r\}, k \in Z_+>$  are linear superspaces 
(generated by indicating in brackets sets of vectors).

Let also $\{\tilde{\phi}_{i,m}\}, \{\tilde{\phi}^{i,m}\}$ are dual bases 
relatively form (\ref{3.5}) bases in superspaces  $\Phi^+, \Phi^-$, respectively. 

We have the following

{\it Proposition 3.3.} The element $R_0$ from proposition 3.1 has the following form 
\begin{equation}
R_0 = exp (\sum_{i,m}(-1)^{deg(\tilde{\phi}_{i,m})}\tilde{\phi}_{i,m}\otimes \tilde{\phi}^{i,m})
\end{equation}

   {\bf Proof}. Let  $B_+ = C[\Phi^+], B_- = C[\Phi^-]$  are commutative function algebras 
on $\Phi^+, \Phi^-$, respectively and 
$\{\tilde{\phi}_{i,m}\}, \{\tilde{\phi}^{i,m}\}$  are above mentioned dual bases. 

Let's fix some total linear ordering of basic and below we'll use notation 
$\{\tilde{\phi}_a \}, \{ \tilde{\phi}^a \}, a \in N.$
Let's prove by induction next formula 
\begin{equation}
<\tilde{\phi}_{i_1}^{n_1} \cdots \tilde{\phi}_{i_k}^{n_k}, (\tilde{\phi}^{i_1})^{m_1} \cdots (\tilde{\phi}^{i_k})^{m_k}> =
\delta_{n_1,m_1} \cdots \delta_{n_k,m_k} n_1! \cdots n_k!
\end{equation}

It is easy to verify the base of induction for 
$k=1, n_1=1$
$<\tilde{\phi}_{i_1}, \tilde{\phi}^{i_1}>=1, <\tilde{\phi}_{i_1}, 1> = 0$.
Further, let
$<\tilde{\phi}_{i_1}^{n}, (\tilde{\phi}^{i_1})^n> = n!$.
Let's show that 
$<\tilde{\phi}_{i_1}^{n+1}, (\tilde{\phi}^{i_1})^{n+1}> = (n+1)!$.
 In fact, \\
$<\tilde{\phi}_i^{n+1}, (\tilde{\phi}^i)^{n+1}> =
<\Delta(\tilde{\phi}_i)\Delta((\tilde{\phi}_i)^n), \tilde{\phi}^i \otimes (\tilde{\phi}_i)^n>= \\
<(\tilde{\phi}_i\otimes 1 + 1\otimes \tilde{\phi}_i)(\sum_{k=0}^{n} C^k_n (\tilde{\phi}_i)^k(\tilde{\phi}_i)^{n-k}), \tilde{\phi}^i \otimes (\tilde{\phi}_i)^n > = \\
(n+1)<\tilde{\phi}_i,\tilde{\phi}^i>
<(\tilde{\phi}_i)^n,(\tilde{\phi}^i)^n > = (n+1)!$.
Using proved formula by induction on $k$ it is proved statement of theorem.
Theorem is proved.

Let now $(f(u))' =\frac{d}{du}(f(u))$. Let's differentiate equality 
(\ref{3.5}) on parameter $u$.  We derive \\
$\frac{d}{du}<\phi^+_i(u), \phi^-_j(u)>= <(\phi^+_i(u))', \phi^-_j(u)> =
\frac{1}{u-v + \frac{1}{2}(\alpha_i, \alpha_j)} - \frac{1}{u-v - \frac{1}{2}(\alpha_i, \alpha_j)}$.
Let $\tilde{\phi}^+_i(u)= \sum_{k=0}^{\infty}\tilde{\phi}_{i,k}u^{-k-1}$,
$\tilde{\phi}^-_i(u)= \sum_{k=0}^{\infty}\tilde{\phi}_{i,-k-1}u^k$. 
Then in terms of generaing functions the pairing 
$<\tilde{\phi}_{i,k}, \tilde{\phi}_{j,l}>=\delta_{ij}\delta_{kl}$
can be rewritten in the next form: 
$<\tilde{\phi}^+_i(u), \tilde{\phi}^-_j(v)> =
\sum_{k,l}<\tilde{\phi}_{i,k}, \tilde{\phi}^{j,l}>u^{-k-1}v^l= \\
(v<1<u)=\delta_{ij}\sum_{k=1}^{\infty}u^{-1}(\frac{v}{u})^k=
\frac{\delta_{ij}}{u-v}$.

Thus we receive that 
\begin{equation}
<\tilde{\phi}^+_i(u), \tilde{\phi}^-_j(v)> =
\frac{\delta_{ij}}{u-v}     \label{3.p}
\end{equation}

Let' introduce a generating vector-functions 

$\tilde{\phi}^{\pm}(u)=
\begin{pmatrix}
\tilde{\phi}^{\pm}_1(u)\\ \tilde{\phi}^{\pm}_2(u)\\ \cdots \\
\tilde{\phi}^{\pm}_r(u)
\end{pmatrix}$

    Then pairing  (\ref{3.p}) we can rewrite in the next matrix equality:
\begin{equation}
<(\tilde{\phi}^+(u))^T, \tilde{\phi}^-(u)> = \frac{E_r}{u-v},
\end{equation}
where $E_r$  be a unit $r\times r$-matrix.

        Let $T: f(v) \rightarrow f(v-1)$  be a shift operator. Clearly that 
\begin{eqnarray*}
<(\phi^-_i(v))', \phi^+_j(v)>=
\frac{1}{u-v + \frac{1}{2}(\alpha_i, \alpha_j)} - \frac{1}{u-v - \frac{1}{2}(\alpha_i, \alpha_j)}=  \\
(id\otimes (T^{b_{ij}}-T^{-b_{ij}}))\frac{\delta_{ij}}{u-v}=
 <\tilde{\phi}^-_i(v), (T^{b_{ij}}-T^{-b_{ij}})\tilde{\phi}_j^+(u)>.
\end{eqnarray*}

Here $b_{ij}=\frac{1}{2}a_{ij}=\frac{1}{2}(\alpha_i,\alpha_j)$.
        Let $A = (a_{ij})_{i,j=1}^r$  be a symmetric Cartan matrix of Lie superalgebra $A(m,n)$.
$\mathfrak{g}$, i. e.  $a_{ij} = (\alpha_i, \alpha_j)$. Let also 
$A(q)=(a_{ij}(q))_{i,j=1}^r$  be a q-analog of Cartan matrix, where
$a_{ij}(q)=[a_{ij}]_q = [(\alpha_i, \alpha_j)]_q=
\frac{q^{a_{ij}}-q^{-a_{ij}}}{q-q^{-1}}.$

Let also $D(q)$ be an inverse matrix to $A(q)$ and $A^T$ denote a transposition of matrix $A$. 
Then we can rewrite previous equality in the next matrix form   
\begin{eqnarray*}
<(\phi^+(v))^T, \phi^-(u)> =
<(\tilde{\phi}^+(u))^T, A(T^{-\frac{1}{2})}(T-T^{-1})\tilde{\phi}^-(v)>
\end{eqnarray*}
Therefore 
\begin{eqnarray}
&<(\tilde{\phi}^-(v))^T, \tilde{\phi}^+(u)>& = \nonumber \\
&<((T-T^{-1})^{-1}D(T^{-\frac{1}{2}})(\phi^-(v))^T, (\phi^+(u))'>&
\end{eqnarray}

   Thus we have the next equality 
\begin{eqnarray}
\frac{E_r}{u-v} =
<((T-T^{-1})^{-1} D(T^{-\frac{1}{2}})(\phi^-(v)), ((\phi^+(u))')^T>
\end{eqnarray}

   So we have diagonalize pairing.   Let's present matrix $D(q)$ in the form 
$D(q)= \frac{1}{[l(\mathfrak{g})]}C(q)$, where $C(q)$ be a matrix with matrix coefficirnts 
being  polynomials of $q$ and $q^{-1}$  with positive integer coefficients 
(i.e. $c_{ij} \in Z[q,q^{-1}]$). Let also  $l(\mathfrak{g})=\check h(\hat{\mathfrak{g}})$  be 
a dual Coxeter number. In these notations the previous formula can be written in the next form

\begin{eqnarray}
&\frac{E_r}{u-v}= &\nonumber \\
&<((T^{l(\mathfrak{g})}-T^{-l(\mathfrak{g})})^{-1}C(T^{-\frac{1}{2}})(\phi^-(v)),
((\phi^+(u))')^T>&
\end{eqnarray}

  From this equality imply formula for the term $R_0$  in the factorizable formula for 
the universal $R$-matrix.

{\bf Theorem 3.1.}
\begin{eqnarray}
&R_0 = &    \nonumber \\
&\prod_{n \ge 0} \exp\sum_{i,j=1}^r \sum_{k \ge 0} ((\phi^+_i(u))')_k \otimes
c_{ji}(T^{-\frac{1}{2}})(\phi^-_j(v+(n+\frac{1}{2})l(\mathfrak{g})))_{-k-1}& \qquad
\end{eqnarray}

\section{Computation of the Universal $R$- matrix of the Yangian $Y(\mathfrak{g})$.}

 First of all let's consider the classical analogs of the argumrnts which will be leaded below. 
Classical $r-$matrix $r(u,v)$ of the classical double \\
$(\mathfrak{g}((u^{-1}))), u^{-1}\mathfrak{g}[[u^{-1}]], \mathfrak{g}[u])$
of the current algebra $\mathfrak{g}[u]$ has next form:
$r(u,v)= \sum_{i,k} e_{i,k} \otimes e^{i,k}$,
where $\{e_{i,k} = e_i u^k\}, \{e^{i,k} = e^i u^{-k-1}\}$ are the dual bases in the 
$\mathfrak{g}[u]), u^{-1}\mathfrak{g}[[u^{-1}]]$, respectively,  with relate to pairing
$$ <f,g>= res(f(u), g(u)),$$
where $(\cdot, \cdot)$ be an invariant bilinear form on $\mathfrak{g}$  and 
$\{e_i\}, \{e^i\}$  are  
dual bases in $\mathfrak{g}$  relative to this form.

It is easy to see that \\
$r = \sum_i \sum_{k=0}^{\infty}e_i\cdot u^k \otimes e^i\cdot
v^{-k-1} = \sum_{k+0}^{\infty} \sum_i e_i \otimes e^i \cdot
v^{-1}(\frac{u}{v})^k = (u<1<v) = \sum_i e_i \otimes e^i
\frac{v^{-1}}{1 - u/v} = \frac{\mathfrak{t}}{v-u}, $\\
where  $\mathfrak{t} = \sum_i e_i \otimes e^i$  be an Casimir operator of universal enveloping 
superalgebra  $U(\mathfrak{g})$  of Lie superalgebra $\mathfrak{g}=
A(m,n)$. Thus we have that 
\begin{equation}
r = \frac{\mathfrak{t}}{v-u}
\end{equation}
Note that this classical  $r$-matrix don't belong to 
$\mathfrak{g}[t]^{\otimes2}$. Let's introduce shift operator 
$T_{\lambda}: f(u) \rightarrow f(u+\lambda)$.  Let's act 
by operator $id \otimes T_{\lambda}$ on $r$. We derive \\
$(id \otimes T_{\lambda})r(u,v) = \frac{\mathfrak{t}}{\lambda -
(u-v)} = \frac{\mathfrak{t}}{\lambda(1 - \lambda^{-1}(u-v))} = \\
\sum_{k=0}^{\infty} \mathfrak{t}(u-v)^k \lambda^{-k-1} =
\sum_{k=0}^{\infty}r_k \lambda^{-k-1}, $ \\
where $r_k \in \mathfrak{g}[t]^{\otimes2}$.
     
   We'll derive this arguments another equivalent way in order that to do the analogy with 
quantum case more  evident. \\
$(id \otimes T_{\lambda})r(u,v)=
\sum_i \sum_{k=0}^{\infty}e_i\cdot u^k \otimes e^i\cdot
(v+\lambda)^{-k-1} = \\
\sum_i \sum_{k=0}^{\infty}e_i\cdot u^k \otimes e^i\cdot
\frac{1}{(\lambda(1 - (-v/\lambda)))^{k+1}} =
\\ \sum_i \sum_{k=0}^{\infty}\sum_{m=0}^{\infty}e_i\cdot u^k \otimes e^i\cdot
(-1)^m C^k_{m+k}v^m\lambda^{-m-k-1}= (n=m+k) =  \\
 \sum_{n=0}^{\infty}\sum_{k=0}^{\infty}\sum_i (-1)^{n-k}C^k_n e_iu^k\otimes e^i
v^{n-k}\lambda^{-n-1}=
\sum_{k=0}^{\infty}\sum_i e_i \otimes e^i(u-v)^n \lambda^{-n-1}$

	Let's try to repeat this argument in the quantum case keeping in mind that Yangian 
be a quantization of the bisuperalgebra of polynomial current $\mathfrak{g}[t]$, 
Yangian Double be a quantization of classical double $\mathfrak{g}((t))$ and universal 
$R$-matrix of Yangian Double be a quantum analog of classical $r$-matrix $r$, and 
above considered $r$-matrix $(id \otimes T_{\lambda})r(u,v)$ is that classical analog of the 
universal $R$- matrix of Yangian wgich we are going to compute. 

     Let's define homomorphism  $T_{\lambda}$ in the quantum case 
$$T_{\lambda}: Y(\mathfrak{g}) \rightarrow Y(\mathfrak{g})$$
by formulas: $T_{\lambda}(x)=x$ for $x \in \mathfrak{g}$,
$T_{\lambda}(a_{i,1}) = a_{i,1} + \lambda a_{i,0}$ for $a \in \{e,f,h\}$.

{\bf Proposition 4.1.}  The action of $T_{\lambda}$ on the generators $a_{i,n}, a \in \{e,f,h\}, n \in Z,$
of Yangian Double $DY(\mathfrak{g})$ is defined by next formulas:
\begin{eqnarray}
T_{\lambda}a_{i,n} = \sum_{k=0}^n C^k_n a_{i,k}\lambda^{n-k}, n \in Z_+  \label{4.6} \\
T_{\lambda}a_{i,-n} = \sum_{k=0}^{\infty}(-1)^k C^{n-1}_{k+n-1} a_{i,k}\lambda^{-n-k-1}, n \in N \label{4.7}
\end{eqnarray}

{\it Proof.}  First of all note that the same formulas define action of   
$T_{\lambda}$ in classical case also.  Let $\tilde{h}_{i,1} = h_{i,1} - \frac{1}{2} h_{i,0}^2$.
 Then it is fulfilled next relations in the Yangian Double $DY(\mathfrak{g})$.
\begin{eqnarray}
&[\tilde{h}_{i,1}, e_{j,n}]=a_{ij}e_{i,n+1}, [\tilde{h}_{i,1}, f_{j,n}]= -a_{ij}f_{i,n+1}, \label{4.8}\\
&[e_{i,k}, f_{j,m}]=\delta_{ij}h_{j, k+m}. \label{4.9} \quad
\end{eqnarray}
Sufficiently to check this formulas for generators $e_{j,n}, f_{j,n}$. 
Because  the relations (\ref{4.7}), (\ref{4.8}) for generators $h_{j,n}$ imply for these relations for 
generators $e_{j,n}, f_{j,n}$ and formula \ref{4.9}. Let's prove the relations (\ref{4.7}), (\ref{4.8}) for generators $e_{j,n}$. For generators $f_{j,n}$ the arguments are the same. Let's prove the formula (\ref{4.7}). For $n=0$ formula (\ref{4.7}) is true on definition. Let this formula is true for all $k \le n$. Let's prove that this formula is true for $k= n+1.$  Let $j$ such that $a_{j,i} \ne 0$. Let's note that  $T_{\lambda}$  be a homomorphism. Then \\
$T_{\lambda}(e_{i,n+1})=T_{\lambda}(a_{ji}^{-1}[\tilde{h}_{j,1}, e_{i,n}])=
a_{ji}^{-1}[T_{\lambda}\tilde{h}_{j,1}, T_{\lambda}e_{j,n+1}]=\\
a_{ji}^{-1}[\tilde{h}_{j,1} + \frac{1}{2}h_{j,0}, \sum_{k=0}^n C^k_n e_{i,k}\lambda^{n-k}]= \\
a_{ji}^{-1}(\sum_{k=0}^n C^k_n[\tilde{h}_{j,1}, e_{i,k}]\lambda^{n-k} +\sum_{k=0}^n C^k_n[h_{j,0}, e_{i,k}]\lambda^{n+1-k})=\\
\sum_{k=0}^{n+1} C^k_{n+1}e_{i,k}\lambda^{n+1-k}.$\\
Formula (\ref{4.7}) is proved.  Let's prove the formula (\ref{4.8}). Let $n=1$. Then
$a_{ji}^{-1}[\tilde{h}_{j,1}, e_{i,-1}]= e_{i,0}$. Let's act on the left-hand  and right-hand sides 
by the operator $T_{\lambda}$. First let's act on the left-hand side. We have \\
$a_{ji}^{-1}[T_{\lambda}(\tilde{h}_{j,1}), T_{\lambda}(e_{i,-1})]=
a_{ji}^{-1}[\tilde{h}_{j,1} + \frac{1}{2}h_{j,0}, \sum_{k=0}^{\infty}(-1)^k C^0_k e_{i,k}\lambda^{-k-2}]=\\
a_{ji}^{-1}[\tilde{h}_{j,1} + \frac{1}{2}h_{j,0}, \sum_{k=0}^{\infty} e_{i,k}\lambda^{-k-2}]=\\
\sum_{k=0}^{\infty}(-1)^k e_{i,k+1}\lambda^{-k-2} + \sum_{k=0}^{\infty}(-1)^k e_{i,k}\lambda^{-k-1} =\\
e_{i,0} + \sum_{k=0}^{\infty}((-1)^k +(-1)^{k+1})e_{i,k+1}\lambda^{-k-2}= e_{i,0}$\\
As the right-hand side don't change under the action $T_{\lambda}$  on definition then we have 
that left-side and right-side hands are coincide. Formula for $n=1$ is verifyed. Let 
formula is proved for all $k \le n$. Let's prove it for $k=n+1$. Let's act as above by operator 
$T_{\lambda}$  on left-hand and right-hand sides of formula 
$a_{ji}^{-1}[\tilde{h}_{j,1}, e_{i,-n-1}]=e_{i,-n}$.  Let's act on the left-hand side. 
$a_{ji}^{-1}[T_{\lambda}\tilde{h}_{j,1}, T_{\lambda}e_{i,-n-1}] =
a_{ji}^{-1}[\tilde{h}_{j,1} + \frac{1}{2}h_{j,0}, \sum_{k=0}^{\infty}(-1)^k C^{n}_{k+n} e_{i,k}\lambda^{-n-k-2}] =
\sum_{k=0}^{\infty}(-1)^k C^{n}_{k+n} e_{i,k+1}\lambda^{-n-k-2} + \sum_{k=0}^{\infty}(-1)^k C^{n}_{k+n} e_{i,k}\lambda^{-n-k-1}=\\
e_{i,0} + \sum_{k=0}^{\infty}(-1)^{k+1} (-C^{n}_{k+n}+C^{n}_{k+n+1})e_{i,k}\lambda^{-n-k-2}= \\
\sum_{k=0}^{\infty}(-1)^k C^{n-1}_{k+n-1} e_{i,k}\lambda^{-n-k-1}= T_{\lambda}e_{i,-n}$\\
We have that left-hand side is equal to right-hand side, and therefore the formula (\ref{4.8}) for all natural numbers $n$ is proved by induction. Proposition is proved. \\

{\it Remark } Note that series defining the value of operator $T_{\lambda}$ on generators 
$a_{i,-n}$ converges for enough large values of $\lambda$.

Now we can to calculate the universal $R$-matrix ${\it R}(\lambda)$ of Yangian 
$Y(\mathfrak{g})$ by formula:
\begin{equation}
{\it R}(\lambda) = (id \otimes T_{-\lambda})R,
\end{equation}
where $R$ be an universal $R$-matrix of double $DY(\mathfrak{g})$. As $R= R_+ R_0 R_-$, then 
acting by operator $id \otimes T_{\lambda}$ on $R$  and using the fact that $T_{\lambda}$ be 
a homomorphism and therefore $id \otimes T_{\lambda}$ be a homomorphism also, we have,
\begin{equation}
{\it R}(\lambda)= {\it R}_+(\lambda) {\it R}_0(\lambda) {\it R}_-(\lambda),
\end{equation}
where ${\it R}_+(\lambda)= (id \otimes T_{-\lambda})R_+, {\it R}_0(\lambda)= (id \otimes T_{-\lambda})R_0,
{\it R}_-(\lambda)= (id \otimes T_{-\lambda})R_-.$

 Note that
\begin{equation}
{\it R}(\lambda)= 1 + \sum_{k=0}^{\infty}{\it R}_k \lambda^{-k-1},
\end{equation}
where $1$ be an unit element in $Y(\mathfrak{g})^{\otimes 2}$,  
${\it R}_k \in Y(\mathfrak{g}) \otimes Y(\mathfrak{g})$.
Such form is not enough suitable as coefficients ${\it R}_k$ have heavy visible form.
Becouse final result let's present in other more visible form. 
 Let's act by operator $id \otimes T_{-\lambda}$ on right-hand side of \ref{3.3}.
We have
\begin{equation}
{\it R}_+(\lambda)  =  \overrightarrow{\prod}_{\beta \in \Xi_+} \exp(-(-1)^{\theta(\beta)} a(\beta) e_{\beta} \otimes T_{-\lambda} e_{-\beta}), \label{4.12}
\end{equation}

Let's calculate separately element $T_{-\lambda} e_{-\beta}$. As $\beta = \beta' + n \delta$, then 
in view \ref{4.7}  we have 
\begin{equation}
{\it R}_-(\lambda) =  \overleftarrow{\prod}_{\beta \in \Xi_-} \exp(-(-1)^{\theta(\beta)} a(\beta) (\sum_{k=0}^{\infty}(-1)^k C^{n-1}_{k+n-1}(e_{\beta}\otimes e_{-\beta +(n+m)\delta})\lambda^{-n-k-1}), \label{4.13}\\
\end{equation}

    Similarly it is calculated element ${\it R_-}(\lambda)$. Summarizing stated  above   
we get 

{\bf Proposition 4.2.}  Terms ${\it R_+}(\lambda),
{\it R_-}(\lambda)$ of universal $R$-matrix of Yangian have the next form 
\begin{eqnarray}
&{\it R}_+(\lambda) =   \overrightarrow{\prod}_{\beta \in \Xi_+} \exp(-(-1)^{\theta(\beta)} a(\beta) (\sum_{k=0}^{\infty}(-1)^k C^{n-1}_{k+n-1}(e_{\beta}\otimes e_{-\beta +(n+m)\delta})\lambda^{-n-k-1}),\qquad \label{4.14}\\
&{\it R}_-(\lambda)  =   \overleftarrow{\prod}_{\beta \in \Xi_-} \exp(-(-1)^{\theta(\beta)} a(\beta) (\sum_{k=0}^{\infty}(-1)^k C^{n-1}_{k+n-1}(e_{\beta}\otimes e_{-\beta +(n+m)\delta})\lambda^{-n-k-1}),\qquad \label{4.15}
\end{eqnarray}

{\it Example 4.1.} Let's consider  example of calculation of terms ${\it R_+}(\lambda),
{\it R_-}(\lambda)$ in the case of the simple Lie algebra $\mathfrak{sl}_2$.
Then\\
${\it R}_+(\lambda)= \overrightarrow{\prod}_{n \ge 0}
\exp(-e_n \otimes T_{\lambda}(f_{-n-1} = \\
\overrightarrow{\prod}_{n \ge 0} \exp((-1)^n \sum_{m=0}^{\infty}
C^n_{m+n} e_n \otimes f_m \lambda^{-n-m-1}) = \\
\exp(\sum_{n=0}^{\infty}  \sum_{m=0}^{\infty}C^n_{m+n}(-1)^n e_n \otimes f_m
\lambda^{-n-m-1})= \exp(\sum_{k=0}^{\infty} (\sum_{m=0}^{k}C^n_{k}(-1)^n e_n \otimes
f_{k-n})\lambda^{-k-1})= \\
\overrightarrow{\prod}_{n \ge 0} \exp((\sum_{m=0}^{k}C^n_{k}(-1)^n e_n \otimes
f_{k-n})\lambda^{-k-1})$ \\
 Similarly it can be calculated the term ${\it R_-}(\lambda)$. Thus we get the next 
formulas 
\begin{eqnarray}
&{\it R}_+(\lambda)= \overrightarrow{\prod}_{n \ge 0} \exp((\sum_{m=0}^{k}C^n_{k}(-1)^n e_n \otimes
f_{k-n})\lambda^{-k-1}), \quad\\
&{\it R}_-(\lambda)= \overleftarrow{\prod}_{n \ge 0} \exp((\sum_{m=0}^{k}C^n_{k}(-1)^n f_n \otimes
e_{k-n})\lambda^{-k-1}), \quad
\end{eqnarray}

    Let  $ord(\beta):= n$, if $\beta = \beta' + n \delta,
\beta' \in \Delta_+(\mathfrak{g})$. Then the proposition 4.2 we can rewrite in the next  form
\begin{eqnarray}
&{\it R}_+(\lambda) =
\overrightarrow{\prod}_{\beta \in \Xi_+} \exp(-(-1)^{\theta(\beta)}
a(\beta) (\sum_{k=0}^{\infty} (-1)^k \quad \nonumber\\
&C^{ord(\beta)-1}_{k + ord(\beta)-1}(e_{\beta} \otimes
e_{-\beta +(ord(\beta)+k)\delta})\lambda^{-ord(\beta)-k-1})),\quad \label{4.16}\\
&{\it R}_-(\lambda) =  \overleftarrow{\prod}_{\beta \in \Xi_-} \exp(-(-1)^{\theta(\beta)} a(\beta)
(\sum_{k=0}^{\infty} \quad \nonumber\\
&(-1)^k C^{ord(\beta)-1}_{k + ord(\beta)-1}(e_{\beta}\otimes e_{-\beta +(ord(\beta)+k)\delta})\lambda^{-ord(\beta)-k-1}), \label{4.17}
\end{eqnarray}

Let's calculate the term ${\it R}_0(\lambda)$. For this it is required next 

{\bf Proposition 4.3.} Shift operator acts on generating function of Cartan generators in the 
following way 
\begin{equation}
T_{\lambda}(h^-_i(u)) = h^+_i(u + \lambda)
\end{equation}

  {\it Proof.} $T_{\lambda}(h^-_i(u)) = T_{\lambda}(1-  \sum_{k=0}^{\infty}h_{i, -k-1}u^k)
= 1 - \sum _{k=0}^{\infty}T_{\lambda}(h_{i, -k-1})u^k =
1 - \sum_{k=0}^{\infty}\sum_{m=0}^{\infty}
C^k_{m+k}h_{i,m})\lambda^{-m-k-1}u^k = 1 + \sum_{m=0}^{\infty}
h_{i,m}(\lambda^+ u)^{-m-1} = h^+_i(u+ \lambda).$ Proposition is proved.\\

{\bf Corollary 4.2.} $T_{\lambda}(\varphi^-_i(u)) = \varphi^+_i(u+ \lambda)$.

{\it Proof.} $T_{\lambda}(\varphi^-_i(u))=
T_{\lambda}(\ln(h^-_i(u)))= \ln(T_{\lambda}(h^-_i(u)))= \ln(h^+_i(u+
\lambda))= \varphi^+_i(u+ \lambda)$.

Now we can calculate the term  ${\it R}_0(\lambda)$.

{\bf Proposition 4.4.} Term ${\it R}_0(\lambda)$ has next form 
\begin{equation}
{\it R}_0(\lambda) =
\prod_{n \ge 0} \exp(\sum_{i,j \in I} \sum_{k \ge 0} ((\phi^+_i(u))')_k \otimes
c_{ji}(T^{-\frac{1}{2}})(\phi^+_j(v+(n+\frac{1}{2})l(\mathfrak{g})+\lambda))_{k})
\label{4.18}
\end{equation}

{\it Example 4.2.} In the case of the simple Lie algebra $\mathfrak{sl}_2$
formula  \ref{4.18} admit the next simple form 
$${\it R}_0(\lambda) =
\prod_{n \ge 0} \exp(-\sum_{k \ge 0}((\phi^+(u))')_k \otimes
(\phi^+(v-2n-1+\lambda))_k $$

   We can present results of this section in the form of the next theorem.

{\bf Theorem 4.4.}   Universal R-matrix of Yangian $Y(A(m,n)$ has the next form 
\begin{equation}
{\it R}(\lambda)= {\it R}_+(\lambda) {\it R}_0(\lambda) {\it R}_-(\lambda),
\end{equation}
where terms ${\it R}_+(\lambda), {\it R}_0(\lambda), {\it R}_-(\lambda)$  is described by, 
respectively,  formulas (\ref{4.16}), (\ref{4.17}), (\ref{4.18}).

\vspace{0.5cm}
    {\bf Acknowledgments.}
        Author are thankful to S.M. Khoroshkin for useful discussions of this work.

\vspace{1cm}

 Math.Dep., Don State Technical University, Gagarin sq.,1, Rostov-na-Donu, 344010, Russia. 

\end{document}